\documentclass[oneside,10pt]{article}          
\usepackage[b5paper]{geometry}	    
\usepackage{amsfonts,amsmath,latexsym,amssymb} 
\usepackage{theorem}                
\usepackage{mathrsfs,upref}         
\usepackage{mathptmx}		    
\usepackage{elemath2024}	            

\theoremstyle{definition}

\newtheorem{definition}{Definition}

\def\Rdst{{\Rst^n}}
\def\cG{\Rdst}
\def\mcM{{\mathcal{M}}}
\newcommand\Xrho{\Xsp_\rho}

\def\Rst{{\mathbb R}}
\def\Rdst{{{\Rst^d}}}
\def\normta#1#2{{  \| {#1}   \|_{#2} \, }}
\def\vol{{\operatorname{vol}}}

\newcommand\ARn{{\Asp(\Rst^n)}}     
\newcommand\Bisp{{\Bsp^1}}     
\newcommand{\Bsp}{{\boldsymbol B}}     
\newcommand\BispN{{(\Bisp, \, \|\ebbes\|^{(1)})}}     
\newcommand{\ebbes}{\mbox{$\,\cdot\,$}}     
\newcommand{\BspN}{(\Bsp, \, \|\ebbes\|_\Bsp)}     
\newcommand\Btsp{{\Bsp^2}}     
\newcommand\BtspN{(\Bsp^2, \, \|\ebbes\|^{(2)})} 
\newcommand{\CORdN}{{\big( \COsp(\Rst^d), \, \|\ebbes\|_\infty \big)}}     
\newcommand{\COsp}{{\Csp_{\negthinspace 0}}}     
\newcommand{\Csp}{{\boldsymbol C}}     
\newcommand{\FF}{{\mathcal{F}}}     
\newcommand\FLi{{\mathcal F}{\negthinspace \Lisp}}     
\newcommand{\Lisp}{{\Lsp^1}}     
\newcommand{\FLiRd}{{ \FLi(\Rdst) }}     
\newcommand\FLiRdN{\big( \FLiRd, \, \|\ebbes\|_{\FLisp} \big)}     
\newcommand{\FLisp}{{ {\mathcal F} \negthinspace \Lisp}}     
\newcommand\FLpsp{{\FF \negthinspace \Lpsp}}     
\newcommand{\Lpsp}{{\Lsp^p}}     
\newcommand{\FT}{{\operatorname{{\mathcal F}}}}     
\newcommand\Hs{{{\cal H}_s}}     
\newcommand\HsRdN{{({\cal H}_s (\Rdst), \, \|\ebbes\|_{\Hs} \big)}}     
\newcommand{\LiRd}{{\Lisp \nth (\Rst^d)}}     
\newcommand\nth{\negthinspace}     
\newcommand{\LiRdN}{\big( \LiRd, \, \|\ebbes\|_1 \big)}     
\newcommand\Liloc{ {\Lsp^1_{ \negthinspace \it{loc}} } }     
\newcommand{\Lsp}{{\boldsymbol L}}     
\newcommand\LilocRd{ \Liloc({\Rdst})}     
\newcommand\LinfRd{ \Lsp^{\infty}(\Rst^{d})}     
\newcommand{\LinfRdN}{\big( \LinfRd, \, \|\ebbes\|_\infty \big)}     
\newcommand{\LiwRd}{{\Liwsp(\Rst^d)}}     
\newcommand{\Liwsp}{{\Lsp^1_{\negthinspace w}}}     
\newcommand{\LiwRdN}{\big( \LiwRd, \, \|\ebbes\|_{1,w} \big)}     
\newcommand{\LpR}{{\Lpsp(\Rst)}}     
\newcommand{\LpRd}{{\Lpsp(\Rst^d)}}     
\newcommand{\LpRdN}{\big( \LpRd, \, \|\ebbes\|_p \big)}     
\newcommand{\LpT}{{\Lpsp(\Tst)}}     
\newcommand{\Tst}{\mathbb{T}}     
\newcommand{\LpTN}{\big( \LpT, \, \|\ebbes\|_p \big)}     
\newcommand{\LtR}{{\Ltsp(\Rst)}}     
\newcommand{\Ltsp}{{\Lsp^2}}     
\newcommand{\LtRd}{{\Ltsp(\Rst^d)}}     
\newcommand{\LtRdN}{\big( \LtRd, \, \|\ebbes\|_2 \big)}     
\newcommand{\Msp}{{\boldsymbol M}}     
\newcommand{\MbRd}{{\Mbsp(\Rst^d)}}     
\newcommand{\Mbsp}{{\Msp_{\negthinspace b}}}     
\newcommand{\SORd}{{\SOsp(\Rst^d)}}     
\newcommand{\SOsp}{{\Ssp_{\negthinspace 0}}}     
\newcommand{\SORdN}{\big( \SORd, \|\ebbes\|_\SOsp \big)}     
\newcommand{\Ssp}{{\boldsymbol S}}     
\newcommand{\ScPRd}{{\ScPsp(\Rst^d)}}     
\newcommand{\ScPsp}{{\Scsp'}}     
\newcommand\ScPRn{{\ScPsp(\Rst^n)}}     
\newcommand{\Scsp}{{\boldsymbol{\mathcal S}}}     
\newcommand{\ScRd}{{\Scsp(\Rst^d)}}     
\newcommand\WCOlpRd{{\Wsp(\COsp,\lpsp)(\Rdst)}}     
\newcommand{\Wsp}{{\boldsymbol W}}     
\newcommand{\lpsp}{{\lsp^p}}     
\newcommand\WFLili{{\Wsp(\FLi \nnth,\lisp)}}     
\newcommand{\nnth}{{ \negthinspace \: \negthinspace }}     
\newcommand{\lisp}{{\lsp^1}}     
\newcommand\WFLiliRd{{\Wsp(\FLi \nth, \lisp)(\Rdst)}}     
\newcommand\WFLiliRdN{{(\WFLiliRd, \|\ebbes\|_{\WFLili})}}     
\newcommand{\WFLplq}{{\Wsp(\FT\Lsp^p,\lqsp)}}     
\newcommand\lqsp{{ \lsp^q}}     
\newcommand\cOsp{{\csp_0}}     
\newcommand{\csp}{{\boldsymbol c}}     
\newcommand{\cOspN}{\big( \cOsp, \|\ebbes\|_{\infty} \big)}     
\newcommand\hatf{{\widehat{f}}}     
\newcommand\hkr{\hookrightarrow}     
\newcommand{\intR}{\int_{\Rst}}     
\newcommand{\intRd}{\int_{\Rst^d}}     
\newcommand{\linsp}{{\lsp^\infty}}     
\newcommand{\lsp}{{\boldsymbol\ell}}     
\newcommand{\linspN}{\big( \linsp, \|\ebbes\|_{\infty} \big)}     
\newcommand\nARn{{\| \cdot | \Asp(\Rst^n)\| }}     
\newcommand\qandq{{ \quad \mbox{and} \quad }}     
\newcommand\suth{{ \, | \, } }     
\def\Xsp{{\boldsymbol X}}
\def\Bnorm#1{{ \| #1 \|_\Bsp }}
\def\by{{\mathbb y}}
\def\NSPB#1{{\left (#1,\| \ebbes \|_{#1} \right )}}
\def\Asp{{\boldsymbol A}}
\def\Rdst{{{\Rst^d}}}
\def\Rst{{\mathbb R}}
\def\cG{\Rdst}

\begin{document}


\title[World of Function Spaces]{Upcoming Terminology in the World of Function Spaces}


\author{Hans G. Feichtinger}

\address{Hans G. Feichtinger, Faculty of Mathematics, University of Vienna,
Oskar-Morgenstern-Platz 1, 1090 Wien \newline
and Acoustic Research Institute, OEAW, AUSTRIA \\
\email{hans.feichtinger@univie.ac.at}
}

\CorrespondingAuthor{Hans G. Feichtinger}


\date{30.09.2024}                               

\keywords{Banach function spaces; BF-spaces; solid; ball function spaces; homogeneous Banach spaces}

\subjclass{Something like: 26D15, 26A51, 32F99, 41A17}



\begin{abstract}
         It is the purpose of this article to compare various concepts of ``function spaces''. In particular we compare notions of the concept
         of {\it Banach Function Spaces} (in the spirit of Luxemburg-Zaanen)
         to the setting of {\it solid BF-spaces} as it is widely used in
         a number of papers in Time-Frequency Analysis or Coorbit Theory.
   More recently ball quasi-Banach function spaces come into the picture.
   In a way this note is a companion to the authors' article entitled
   ``Choosing Function Spaces in Harmonic Analysis'' (from 2015, February
   Fourier Talks at the Norbert Wiener Center, Vol.4) and the recent article
   by E.~Lorist and Z.Nieraeth, entitled ``Banach Function Spaces done right'' (Indag. Math., 2023). We restrict our attention to the Banach
   space case, avoiding the complications arising in the quasi-Banach case.
\end{abstract}

\maketitle

%




\section{Introduction}

Function spaces play a crucial role in the work of Lars-Erik Persson. Among his contributions, he co-founded the Journal of Function Spaces and Applications, now known as the Journal of Function Spaces. He was its Chief Editor from 2003 to 2011.
Given the vast number of function spaces, even when we limit our focus to Banach spaces relevant to Fourier and Time-Frequency Analysis, there is a clear need for guidance on selecting the appropriate space for specific tasks.

This article serves as a complement to \cite{fe15}, which provides an overview of various function spaces discussing a variety of construction principles, and indicating some of the relationships between such families. Its primary goal is to offer a supplementary guide that introduces new terms and concepts that have gained prominence in the last decade and were not included in \cite{fe15}.
Unlike the interesting article \cite{loni23} we avoid the quasi-Banach case.
Partially this is due to a personal choice, but also based on the observation that the extra technical efforts required does not pay back when it comes to the use of such spaces, with very few exceptions.
We definitely share the view of the authors of \cite{loni23} that the Fatou property should not be part of the definition, and that it would be good to return to the original setting of Banach Function Spaces as developed already decades ago by Luxemburg-Zaanen, see \cite{za67}.


\section{Thoughts inspired by Hans Triebel}

If one looks at the publications of H.~Triebel on Function Spaces, in particular \cite{tr83}, \cite{tr92},\cite{tr01}, \cite{tr16} or \cite{tr20},
it is clear that in all his work he talks about Banach spaces of
tempered distributions, in the sense that the Banach space $\BspN$ is
continuously embedded into $\ScPRd$, the space of tempered distributions, or equivalently:
\begin{equation}\label{BembSCP}
  \lim_{n \to \infty} \Bnorm{f_n} = 0 \,\,  \Rightarrow  \,\, \lim_{n \to \infty} f_n(k) = 0
  \quad \forall k \in \ScRd.
\end{equation}
Obviously, classical function spaces such as the space $\LpRdN$, with
$1 \leq p \leq \infty$, or Sobolev spaces $\HsRdN$, for $s \in \Rst$ are
included, but this setting has the big advantage that for any such
Banach space also the dual space belongs to this family, if $\ScRd$
is dense in $\BspN$.


Reflecting on the different norms used in the construction of new
spaces Triebel comes up with the following useful terminology
(see \cite{tr16}, section 1.6, p.5) given below. Recall that
we restrict our attention to the Banach space case:
%
%
\begin{enumerate}
\item A norm\footnote{We use here Triebel's notation.}
  $\nARn$ is called {\bf admissible} (with respect to $\ScPRd$) if it makes sense to test any $f \in \ScPRn $ for whether it belongs to the corresponding space $\ARn$ or not. The elements of $\ARn$ are those with finite norm.
\item A norm  $\nARn$ is called {\bf regional} if it makes sense to test any element of some region (linear subset of $\ScPRn$) for whether it belongs to the corresponding space $\ARn$ or not.
\item Within a given fixed Banach space $\ARn$,
equivalent norms are  called {\bf domestic}. They are not necessarily defining (admissible or regional).
\end{enumerate}

Good examples are the following ones (for illustration):
\begin{itemize}
  \item We start with the $\Lpsp$-norms, with $1 \leq p \leq \infty$. They are {\it regional} on the  subset of regular tempered distributions which are given by locally  integrable functions $f$. One could also take the space  $\mathscr{M}(\mathbb{R}^d)$ of all (equivalence classes of) measurable functions and come up with the same spaces. Note that outside of this region, which allows to apply the Lebesgue integral on a given function, because expressions like  $|f|^p$ (hence an integral over such an object) would not even be meaningful.

      The theory of Banach function spaces can be described as an important class of regional Banach spaces within the realm of (equivalence classes of) measurable functions over some measure space, using
      typically some function space norms $\rho$, satisfying (among others
      $\rho(f) = \rho(|f|)$ and setting
      $$\Xrho := \{ f \in \mcM \, | \, \normta {f} {\Xrho} :=  \rho(f) < \infty \}. $$

      \item The Fourier characterization of Besov spaces is a good example
      of an {\it admissible} norm, starting from the overall reservoir
      $\ScPRd$ of tempered distributions on $\Rdst$. Filtering an arbitrary tempered distributions by dyadic band-pass filters results in continuous functions, and consequently the $\Lpsp$-norms make sense. The corresponding sequence space norms applied to such a sequence
      of $\Lpsp$-norms is finite or infinity, thus defining the
      elements of the (inhomogeneous) Besov space as claimed. There is
      a never ending stream of publications working in this setting.

  \item A good example of a {\it domestic} norms are Wiener amalgam norms
  corresponding to the spaces $\WCOlpRd$, on spaces of band-limited
  functions (or spline-type spaces, say cubic splines in $\LpRd$)
  which are equivalent on those subspaces of $\LpRd$.
  A more naive example would be $\CORdN$ (with the $\sup$-norm)
  which is a closed subspace of $\LinfRdN$, but of course finiteness
  of the $\sup$-norm does not guarantee that a bounded function
  belongs vanishes at infinity.
\end{itemize}

\section{BK-spaces}

Obviously it is easier to deal with function spaces over discrete domains,
which are usually just called {\it sequence spaces}.
The corresponding theory of BK-spaces, i.e. (partially German terminology!)
of ``Banach-Koordinaten R\"aume'' goes probably back to G.~K\"othe. Thus they are also called ``K\"othe sequence spaces''.  These spaces have been studied
in the 60th and 70th and are Banach spaces of sequences with the (natural,
but non-trivial) assumption that coordinate mappings are continuous). In other
words, a Banach space $\BspN$ of sequences over a (discrete) domain $I$
can be described as families $X = (x_i)_{i \in I}$. Thus we assume that
for every $i \in I$ there exists $C_i < \infty$ such that %
$$   |x_i| \leq C_i \Bnorm{ (x_i)_{i \in I}}, \quad \forall  X \in \Bsp. $$

There is an abundance of such BK-spaces. Many of them are {\it solid},
i.e satisfy the property that for any $X \in \Bsp$ and any family
$ Y = (y_i)_{i \in I} $ one has:
$$  |y_i| \leq |x_i| \, \,\,  \forall i \in I  \,\,\,  \Rightarrow \, \,\, Y \in \Bsp  \qandq \Bnorm{Y} \leq \Bnorm{X}.$$
Clearly weighted $\lpsp$-spaces belong to the class of solid spaces,
while the Banach spaces  $\NSPB \FLpsp$ of all Fourier coefficients of functions in
$\LpTN$, with $p \neq 2$, with the norm
$ \normta {\hatf} {\FLpsp} := \normta {f} {\LpT}$
are good examples of BK-spaces which are {\it not solid}. As V.~Temlyakov has shown in \cite{te11}  this has serious consequences for greedy algorithms (the attempt to recover $f$ from the big Fourier coefficients via hard thresholding).

A good example of usefulness of abstract functional analytic arguments
being useful in such a context is the following simple consequence
of the Closed Graph Theorem: Assume that we have two BK-spaces
$\BispN$ and $\BtspN$, such that $\Bisp \subseteq \Btsp$. Then
the embedding is automatically continuous, i.e. there exists
$C > 0$ such that $ \normta f \Btsp \leq C \normta f \Bisp \, , \,\, f \in \Bisp.$

Let us just give a typical list of references
\cite{go61-1},
\cite{sa66},
\cite{sa64-1} and papers of D.J.H.~Garling and others, such as \cite{ga69-2} or \cite{ga66-3}.


\section{Function Spaces in Harmonic Analysis}

Traditionally Fourier Analysis has been the source for many developments
in the early development of Functional Analysis. In the second half of
the 20th century (and up to now) it appears to be clear that the Fourier
transform, viewed as an {\it integral transform}, but also the  important
notion of {\it convolution} (if defined pointwise a.e.) requires
to work with the Lebesgue space $\LiRd$. In fact, $\LiRdN$ is a Banach algebra
with respect to convolution, with bounded approximate units (so-called
Dirac sequences). Plancherel's Theorem on the other hand shows that
(suitably normalized) it can be viewed as an isometric automorphism
on the Hilbert space $\LtRdN$. Consequently, the scale of space
$\LpRdN$, with $1 \leq p \leq \infty$ play an important role in
Fourier Analysis, and many results valid for these spaces have been
extended to more general Banach spaces of (measurable) functions
on $\Rdst$ (or more generally on LCA groups).
As an illustration let me provide a short list of books that have
influenced this author the most:
\begin{enumerate}
  \item H.~Reiter's book \cite{re68} (new edition with I.~Stegeman \cite{rest00});
  \item Y.~Katznelsons's book on Harmonic Analysis \cite{ka68} (last reprint 2004,  \cite{ka04-1});
  \item P.~Butzer and R.~Nessel book \cite{bune71} 
  emphasizes the approximation theoretical side of Fourier Analysis,
      describing the various summability methods;
  \item The  book \cite{st70} by E.~Stein  lays the foundations
   for many modern developments;
  \item The books of J.~Peetre (\cite{pe76}) and H.~Triebel (see above)
     describe the world of Besov-Triebel-Lizorkin spaces from the
     point of view of interpolation theory;
  \item Standard books on interpolation theory at that time have been
     the books of Bergh/L\"ofstr\"om (based on courses by J.~Peetre),
     as well as  Bennett/Sharpley (\cite{besh88});
  \item   G.~Folland's {\it Harmonic Analysis on Phase Space} \cite{fo89}
   already opens the door to what is nowadays called Time-Frequency Analysis.

\end{enumerate}

Reiter's book describes specifically two families of Banach spaces of
functions (continuously embedded into $\LiRdN$),
the weighted $\Lisp$-algebras $\LiwRdN$,
called {\it Beurling algebras} with bounded approximate units, and on the
other hand the so-called {\it Segal algebras} $\BspN$, which are {\it Banach ideals}
in $\LiRdN$, satisfying
\begin{equation}\label{LiIdeal00}
   \Bnorm{g \ast f} \leq \normta {g} {\Lisp}  \Bnorm{f}, \quad
       g \in \Lisp, f \in \Bsp.
\end{equation}
This property is obtained
by viewing $g \ast f$ as vector-valued
$  g \ast f = \intRd  T_y f g(y) dy,$
using the following two properties
\begin{equation}\label{homBan00}
   \Bnorm{T_x f}  = \Bnorm {f} \qandq  \lim_{x \to 0} \Bnorm{T_x f -f} = 0
     \quad \forall  f \in \Bsp.
\end{equation}
It also applies to the more general situation of Banach spaces
$\BspN \hkr \LilocRd$ satisfying (\ref{homBan00}), which are known
as {\it homogeneous Banach spaces} in the book of Y.~Katznelson \cite{ka04-1}.
The article \cite{fe22} resumes this topic and derives the inclusion
(with corresponding norm estimate)   $\MbRd \ast \Bsp \subseteq \Bsp$, based
on the novel approach to convolution provided in \cite{fe17}.

Another set of function spaces, the so-called {\it modulation
spaces}, introduced already in 1983 in \cite{fe83-4} and published
in \cite{fe03}.
The prototypical example $\SORdN$ (introduced in \cite{fe81-2})
is meanwhile known (cf.\ \cite{rest00}) as ``Feichtinger's algebra''.
It is both a Wiener amalgam space (\cite{fe83},  \cite{he03}
or \cite{fe24-3} for a very recent survey), but also a modulation
space. These modulation spaces are Wiener amalgam space of the
form $\WFLplq$ on the Fourier transform side.

Nowadays one can say that {\it modulation spaces}
play an important role in Time-Frequency Analysis,
Gabor Analysis, but also in the modern theory of pseudo-differential
operators. Relevant books are:
 \begin{enumerate}
   \item The book of K.~Gr\"ochenig on the ``Foundations of Time-Frequency
     Analysis'' has become a standard reference (\cite{gr01});
   \item 
The book \cite{beok20} by A.~Benyi and K.~Okoudjou
 connects modulation spaces with   pseudo-differential operators;
   \item The book \cite{coro20} by E.~Cordero and L.~Rodino
    is entitled ``Time-Frequency Analysis of Operators and Applications'';
   \item Several books by M.~de Gosson make use of modulation spaces, in
    particular let us mention
    \cite{de17}, dealing with ``The Wigner Transform''.
    \item There are various sources meanwhile covering the
    connection between function spaces and the metaplectic
    group (starting with \cite{we64}, \cite{fo89}, \cite{re89}).
  \end{enumerate}



\section{Banach Function Spaces}

Moving from the discrete setting of BK-spaces to the continuous setting we observe at the outset that - despite some obvious vague ideas -
the concept of what it means to deal with a {\it  Banach space of functions}
is by no means mathematically well defined and depends on the context and the authors dealing with questions in this area.

While it is clear that  point evaluations $\delta_x:  f \mapsto f(x)$ are
meaningful on spaces of continuous functions they are not anymore meaningful
for standard examples such as $\LpT$ or even for $\LtR$, except for subspaces
of band-limited functions, which form a RKH (reproducing kernel Hilbert space.

The theory of what has been called {\bf Banach Function Spaces} seems to
fit here. It has been
developed in great detail by Luxemburg and Zaanen, in a series of more
than 17 papers, starting in 1963. It was based on the PhD thesis of
Luxemburg  \cite{lu55}. The book \cite{za67}, Chap.15 is a standard
reference here. There are even at least two more books on the subject:
\cite{luza71-1} and
\cite{za83-1}, dealing with so-called {\it Riesz spaces}.  


The goal of this theory of Banach Function Spaces was to have common setting to deal with $\Lpsp$-spaces and their generalizations, such as Orlicz and Lorentz
spaces (see \cite{besh88}) of measurable functions over general measure spaces.
Thus it is natural to include in the definition also the solidity
condition, and also that the indicator function of any set of finite
measure should belong to $\Bsp$. In fact, for rearrangement
invariant spaces this is no problem, but weighted function space
do not contain indicator functions of sets of finite measure, if they
are spread out too much.

Banach function spaces are slightly different from solid BF-spaces.
On the one hand, Banach function spaces are usually assumed to contain
the indicator function of any measurable set of finite measure. Probably
this property is inspired by the fact that one was mostly concentrating
on so-called rearrangement invariant function spaces in the early period.
However, in the context of Banach spaces of measurable functions on locally
compact groups it is more natural to assume on the one hand local integrability
of the functions under consideration (with respect to the Haar measure, left or
right) and the fact that the indicator functions of arbitrary compact subsets
should belong to such a space. This mostly excludes (large) sets of the
underlying group where all the functions from a given Banach space $\BspN$
vanish. In a series of recent papers this condition is prominently appearing
under the name of {\bf ball Banach function spaces} (BBFS).

Unfortunately the concept of the BBFS concept requires (for convenience) the
validity of Fatou's Lemma, which appears to be somehow related to the admissibility condition of Triebel. Anyway, given any Banach function
space its $\alpha$-dual, also called Koethe-dual, namely
$$\Bsp^\alpha := \{ h \in \Liloc \suth hf \in \LiRd \forall f \in \Bsp\}.$$
Any $\alpha$-dual space satisfies the Fatou property (and conversely), but
some important Banach function spaces (in the original sense) do not satisfy
this property. In fact, a few non-reflexive Banach spaces may be excluded in this case. In the discrete setting $\cOspN$ would provide an example. It is
a closed subspace of $\linspN$, which obviously does not satisfy the Fatou property. The authors of \cite{loni23} share this opinion.


\section{BF-spaces in Harmonic Analysis}


For applications in Harmonic Analysis, especially in the study of
translation invariant Banach spaces of functions it was natural
to assume that one has local integrability (with respect to the
Lebesgue measure, which is the Haar measure on the LCA group $\Rdst$).
We thus came up (already in \cite{fe79}) with the concept of BF-spaces,
thinking of ``{\it Banach Function}'' spaces (instead of Banach Koordinaten
spaces). This setting is meanwhile widely used in papers in the field
of time-frequency analysis. It would go beyond the scope of this note
to list only a representative list of papers using this notion.

\begin{definition} \label{BFdef00}
A BF-space  is a Banach space of locally integrable functions over a locally compact group. More precisely,
$\BspN$ is a Banach space, continuously embedded into  $\Liloc(\cG)$, which means that for every compact
subset $ K \subseteq \cG$ there exists a constant $C_K > 0$ such that
\begin{equation}  \label{Lilocemb}
    \int_K  |f(x)| dx  \leq C_K \|f\|_\Bsp   \quad \forall f \in \Bsp.
\end{equation}
\end{definition}

For discrete groups it is sufficient to have the estimate for single points, because compact sets are then
just finite (unions of one point) sets. The corresponding (older) terminology is thus a special case of
the concept of BF-space for discrete groups:

\begin{definition} \label{BFsolid} 
A BF-{\it space}  
is called {\it solid} if the following is true: \newline
Given  $g  \in \Bsp$ and 
$f  \in {\mathscr{M}(\mathbb{R}^d)} $ with $ |f(x)| \leq  |g(x)| $ a.e.
implies
\begin{equation} \label{solidest1}
 f \in \Bsp  \quad \mbox{and}  \quad  \|f\|_\Bsp  \leq  \|g\|_\Bsp. \end{equation}
\end{definition}
An important subclass of solid BF-spaces are the translation invariant ones.

\section{Ball Banach Function Spaces}



Next, we recall the concept of
ball quasi-Banach function spaces
\cite[Definition 2.2]{hosayaya17}, where the motivation for
the change of assumptions compared to the classical Banach
function spaces is given.
They have been used extensively in a series of papers, such as
\cite{yayuzh23,dagrpayayuzh24,tayayuzh21,ho21-1,hosayaya17}.
It is argued that especially
the assumption of local integrability may be too restrictive,
e.g.\ when one wants to include the spaces $\LpRd$ with $0 < p < 1$
in the consideration (family of quasi-Banach spaces). Again, for
simplification we adapt the definition to the Banach space case:
\begin{definition}\label{debqbf}
A normed linear space  $\Xsp \subset \mathscr{M}(\mathbb{R}^d)$,
equipped with a norm $\|\cdot\|_{\Xsp}$ which makes sense
for all functions in $\mathscr{M}(\mathbb{R}^d) $, is called
a \emph{ball Banach function space}
(for short, BBFS)  if
it satisfies the following conditions:
\begin{enumerate}
\item [(i)]
For any $f\in  \mathscr{M}(\mathbb{R}^d)$, 
$\|f\|_{\Xsp}=0$ if and only if $f(x)=0$ almost everywhere.
\item [(ii)] For any $g \in \Xsp$ one has:
Given $f\in \mathscr{M}(\mathbb{R}^d)$ satisfying $|f(x)|\leq |g(x)|$ almost everywhere, one
has $g \in \Xsp$ and $\|g\|_{\Xsp}\leq \|f\|_{\Xsp}$ ({\it solidity}).
\item [(iii)]
For any sequence
$\{f_n\}_{n \in\mathbb{N}}$ in $\Xsp$ 
with $f_n(x) \uparrow f(x)$ a.e.\,
for some $f\in \mathscr{M}(\mathbb{R}^d) $ one has $f \in \Xsp$ and
$\|f_n\|_{\Xsp}\uparrow \|f\|_{\Xsp}$
as $n\to\infty$ ({\it Fatou property});
\item [(iv)]
For any ball $B\subset \Rdst, \,  {\bf{1}}_B \in \Xsp$.
\end{enumerate}
\end{definition}

We believe, that the assumption (\ref{Lilocemb}) is not really restrictive
for applications in Harmonic Analysis, among others. In fact, the
space $\LpRd$ with $ 0 < p < 1$ cannot even be interpreted as (quasi)-
Banach spaces of tempered distributions, since a singularity of
$h \in \LpR$ may prohibit the interpretation as a regular distributions via
$ \sigma_h(k)  := \intR f(x) \, h(x) dx $,
just choose a non-negative bump function $h \in \Scsp(\Rst)$.

\section{Double Modules}

A class of function spaces (of whatever family you may start) which is gaining relevance in the last years is the collection of Banach spaces having a double module structure, namely with respect to both convolution and pointwise multiplication. A systematic study of such spaces was already undertaken
in \cite{brfe83}. Note that the term ``{\it double module}'' has to be distinguished from the term ``bi-module'', where one typically has a left
and a right action of two different Banach algebras, and these two actions
commute, which in such a case typically corresponds to some kind of {\it
associativity law}.

For this author the first encounter with such structures was the study
of {\it strongly character invariant Segal algebra}. This are Segal
algebras in Reiter's sense (hence homogeneous Banach spaces and thus
furthermore Banach modules over $\LiRdN$ with respect to convolution,
in fact Banach ideals) with the extra property that multiplication
with characters is isometric, i.e.\ , $\Bnorm {\chi_s f}  = \Bnorm{f}$
for all $s \in \Rdst$. In this situation one shows that $s \mapsto \chi_s f$
is a continuous mapping and via integration (see \cite{fe22} for a modern
approach) one obtains that $\FLiRd \cdot \Bsp = \Bsp$ defines an additional
Banach module structure.

It was one of the key results of \cite{fe81-2} to demonstrate that there
is a minimal space in this family of Banach spaces (it is no problem to
consider only subspaces of $\LtRdN$, the most important space in this family)
which can be characterized as the Wiener amalgam space $\WFLiliRdN$,
denoted by $\SORdN$ in this paper. Any solid BF-space is a pointwise
Banach module over $\LinfRdN$, hence also over  $\CORdN$ and in particular
over $\FLiRdN$. In fact, in most cases such a double module structure
arises via the {\it integrated group representation} of either the
group $\Rdst$ acting by translations or (the frequency domain) $\Rdst$,
acting via multiplication with the characters $\chi_s$.
One can combine the two actions to the action of the (reduced)
Heisenberg group,
\cite{auen20}, \cite{jalu19}, \cite{jalu18}, \cite{jalu21},
or the very recent article \cite{falapa24}.

By relaxing the condition of isometric action and replaces it by the
assumption that the given Banach space (of tempered distributions)
is simply invariant under the commutative action of translation
operators respectively the modulation operators (multiplication
with characters). The operator norm on $\BspN$ is then of
course a submultiplicative function (a so-called {\it Beurling weight},
hence one may assume that there exists a submultiplicative function $w$ satisfying
$ w(x+y) \leq w(x) w(y)$ for $x,y \in \Rdst$, and with
$$   \Bnorm{T_x f} \leq w(x) \Bnorm{f}, \quad f \in \Bsp, x \in \Rdst. $$
In a similar way we can assume that there is some Beurling weights
such that
$$   \Bnorm{\chi_s f} \leq v(x) \Bnorm{f}, \quad f \in \Bsp, s \in \Rdst. $$
Such space appear in the recent literature under the name of
 {\it translation and modulation invariant function spaces}.
At first one can work with Banach spaces of tempered distributions
and allow weight functions $w(x)$ and $v(s)$ of polynomial growth,
see \cite{diprvi15}, or later extend it to the setting of
ultra-distributions (\cite{dipiprvi16}, \cite{dipiprvi19}, \cite{fegu21},
\cite{fegu22}). The validity of the {\it bounded approximation property}
for such spaces was derived in \cite{fe22-1}.

 \bibliographystyle{abbrv}  

\begin{thebibliography}{30}

\bibitem{auen20}
A.~{A}ustad and U.~{E}nstad.
\newblock {H}eisenberg modules as function spaces.
\newblock {\em J. Fourier Anal. Appl.}, 26(2):24, 2020.

\bibitem{besh88}
C.~{B}ennett and R.~C. {S}harpley.
\newblock {\em {I}nterpolation of {O}perators.}
\newblock {P}ure and {A}pplied {M}athematics, {V}ol. 129. {B}oston etc.:
  {A}cademic {P}ress, 1988.

\bibitem{beok20}
A.~{B}enyi and K.~A. {O}koudjou.
\newblock {\em {M}odulation {S}paces. {W}ith {A}pplications to
  {P}seudodifferential {O}perators and {N}onlinear {S}chr{\"o}dinger
  {E}quations}.
\newblock {A}ppl. {N}um. {H}arm. {A}nal.,
  ({B}irkh{\"a}user), {N}ew {Y}ork, 2020.

\bibitem{brfe83}
W.~{B}raun and H.~G. {F}eichtinger.
\newblock {B}anach spaces of distributions having two module structures.
\newblock {\em J. Funct. Anal.}, 51:174--212, 1983.

\bibitem{bune71}
P.~L. {B}utzer and R.~J. {N}essel.
\newblock {\em {F}ourier {A}nalysis and {A}pproximation. {V}ol. 1:
  {O}ne-dimensional {T}heory}.
\newblock {B}irkh{\"a}user, {S}tuttgart, 1971.

\bibitem{coro20}
E.~{C}ordero and L.~{R}odino.
\newblock {\em {T}ime-{F}requency {A}nalysis of {O}perators and
  {A}pplications}.
\newblock {D}e {G}ruyter {S}tudies in {M}athematics, {B}erlin, 2020.

\bibitem{dagrpayayuzh24}
F.~{D}ai, L.~{G}rafakos, Z.~{P}an, D.~{Y}ang, W.~{Y}uan, and Y.~{Z}hang.
\newblock {T}he {B}ourgain--{B}rezis--{M}ironescu formula on ball {B}anach
  function spaces.
\newblock {\em {M}athematische {A}nnalen}, 388(2):1691--1768, 2024.

\bibitem{de17}
M.~A. de~{G}osson.
\newblock {\em {T}he {W}igner {T}ransform}.
\newblock {A}dvanced {T}extbooks in {M}athematics. {W}orld {S}cientific
  {P}ublishing {C}o. {P}te. {L}td., {H}ackensack, 2017.

\bibitem{dipiprvi16}
P.~{D}imovski, S.~{P}ilipovic, B.~{P}rangoski, and J.~{V}indas.
\newblock {C}onvolution of ultradistributions and ultradistribution spaces
  associated to translation-invariant {B}anach spaces.
\newblock {\em {K}yoto {J}. {M}ath.}, 56(2):401--440, 2016.

\bibitem{dipiprvi19}
P.~{D}imovski, S.~{P}ilipovic, B.~{P}rangoski, and J.~{V}indas.
\newblock {T}ranslation-modulation invariant {B}anach spaces of
  ultradistributions.
\newblock {\em J. Fourier Anal. Appl.}, 25(3):819--841, 2019.

\bibitem{diprvi15}
P.~{D}imovski, B.~{P}rangoski, and J.~{V}indas.
\newblock {O}n a class of translation-invariant spaces of quasianalytic
  ultradistributions.
\newblock {\em Novi Sad J. Math.}, 45(1):143--175, 2015.

\bibitem{falapa24}
C.~{F}arsi, F.~{L}atremoliere, and J.~{P}acker.
\newblock {C}onvergence of inductive sequences of spectral triples for the
  spectral propinquity.
\newblock {\em {A}dvances in {M}athematics}, 437:109442, 2024.

\bibitem{fe79}
H.~G. {F}eichtinger.
\newblock {G}ewichtsfunktionen auf lokalkompakten {G}ruppen.
\newblock {\em {S}itzber. d. {\"o}sterr. {A}kad. {W}iss.}, 188:451--471, 1979.

\bibitem{fe81-2}
H.~G. {F}eichtinger.
\newblock {O}n a new {S}egal algebra.
\newblock {\em Monatsh. Math.}, 92:269--289, 1981.

\bibitem{fe83}
H.~G. {F}eichtinger.
\newblock {B}anach convolution algebras of {W}iener type.
\newblock In {\em {P}roc. {C}onf. on {F}unctions, {S}eries, {O}perators,
  {B}udapest 1980}, volume~35 of {\em {C}olloq. {M}ath. {S}oc. {J}anos
  {B}olyai}, pages 509--524. {N}orth-{H}olland, {A}msterdam, {E}ds. {B}.
  {S}z.-{N}agy and {J}. {S}zabados. edition, 1983.

\bibitem{fe83-4}
H.~G. {F}eichtinger.
\newblock {M}odulation spaces on locally compact {A}belian groups.
\newblock Technical report, {U}niversity of {V}ienna, {J}anuary 1983.

\bibitem{fe03}
H.~G. {F}eichtinger.
\newblock {M}odulation spaces on locally compact {A}belian groups.
\newblock In R.~{R}adha, M.~{K}rishna, and S.~{T}hangavelu, editors, {\em
  {P}roc. {I}nternat. {C}onf. on {W}avelets and {A}pplications}, pages 1--56,
  {C}hennai, {J}anuary 2002, 2003. {N}ew {D}elhi {A}llied {P}ublishers.

\bibitem{fe15}
H.~G. {F}eichtinger.
\newblock {\em {C}hoosing {F}unction {S}paces in {H}armonic {A}nalysis},
  volume~4 of {\em {T}he {F}ebruary {F}ourier {T}alks at the {N}orbert {W}iener
  {C}enter, {A}ppl. {N}umer. {H}armon. {A}nal.}, pages 65--101.
\newblock {B}irkh{\"a}user/{S}pringer, {C}ham, 2015.

\bibitem{fe17}
H.~G. {F}eichtinger.
\newblock {A} novel mathematical approach to the theory of translation
  invariant linear systems.
\newblock In I.~{P}esenson, Q.~{L}e {G}ia, A.~{M}ayeli, H.~{M}haskar, and
  D.~{Z}hou, editors, {\em {R}ecent {A}pplications of {H}armonic {A}nalysis to
  {F}unction {S}paces, {D}ifferential {E}quations, and {D}ata {S}cience.},
  {A}pplied and {N}umerical {H}armonic {A}nalysis., pages 483--516.
  {B}irkh{\"a}user, {C}ham, 2017.

\bibitem{fe22}
H.~G. {F}eichtinger.
\newblock {H}omogeneous {B}anach spaces as {B}anach convolution modules over
  ${M} ({G})$.
\newblock {\em {M}athematics}, 10(3):1--22, 2022.

\bibitem{fe22-1}
H.~G. {F}eichtinger.
\newblock {T}ranslation and modulation invariant {B}anach spaces of tempered
  distributions satisfy the {M}etric {A}pproximation {P}roperty.
\newblock {\em {A}ppl. {A}nalysis}, 20(6):1271--1293, {N}ov 2022.

\bibitem{fe24-3}
H.~G. {F}eichtinger.
\newblock {T}he concept of {W}iener amalgam spaces.
\newblock In {\em {A}kram {A}ldroubi 65th {B}irthday {V}olume}, {A}ppl. {N}um.
  {H}arm. {A}nalysis, pages 1--25. {B}irkh{\"a}user/{S}pringer, 2024.

\bibitem{fegu21}
H.~G. {F}eichtinger and A.~{G}umber.
\newblock {C}ompleteness of shifted dilates in invariant {B}anach spaces of
  tempered distributions.
\newblock {\em Proc. Amer. Math. Soc.}, 149(12):5195--5210, 2021.

\bibitem{fo89}
G.~B. {F}olland.
\newblock {\em {H}armonic {A}nalysis in {P}hase {S}pace}.
\newblock {P}rinceton {U}niversity {P}ress, {P}rinceton, {N}.{J}., 1989.

\bibitem{ga66-3}
D.~J.~H. {G}arling.
\newblock {O}n symmetric sequence spaces.
\newblock {\em Proc. Lond. Math. Soc. (3)}, 16:85--106, 1966.

\bibitem{ga69-2}
D.~J.~H. {G}arling.
\newblock {A} class of reflexive symmetric {B}{K}-spaces.
\newblock {\em Canad. J. Math.}, 21:602--608, 1969.

\bibitem{go61-1}
G.~{G}oes.
\newblock {C}omplementary spaces of {F}ourier coefficients, convolutions, and
  generalized matrix transformations and operators between {B}{K}-spaces.
\newblock {\em {J}. {M}ath. {M}ech.}, 10:135--157, 1961.

\bibitem{gr01}
K.~{G}r{\"o}chenig.
\newblock {\em {F}oundations of {T}ime-{F}requency {A}nalysis}.
\newblock {A}ppl. {N}umer. {H}armon. {A}nal. {B}irkh{\"a}user, {B}oston,
  {M}{A}, 2001.

\bibitem{he03}
C.~{H}eil.
\newblock {A}n introduction to weighted {W}iener amalgams.
\newblock In M.~{K}rishna, R.~{R}adha, and S.~{T}hangavelu, editors, {\em
  {W}avelets and their {A}pplications ({C}hennai, {J}anuary 2002)}, pages
  183--216. {A}llied {P}ublishers, {N}ew {D}elhi, 2003.

\bibitem{ho21-1}
K.-P. {H}o.
\newblock {E}rdelyi--{K}ober fractional integral operators on ball {B}anach
  function spaces.
\newblock {\em {R}endiconti {S}eminario mat. {U}niv. di
  {P}adova}, 145:93--106, 2021.


\bibitem{jalu19}
M.~{J}akobsen and F.~{L}uef.
\newblock {S}ampling and periodization of generators of {H}eisenberg modules.
\newblock {\em {I}nternational {J}ournal of {M}athematics}, 30(10):1950051,
  2019.

\bibitem{jalu18}
M.~S. {J}akobsen and F.~{L}uef.
\newblock {D}uality of {G}abor frames and {H}eisenberg modules.
\newblock {\em {A}r{X}iv}:1806.05616, June 2018.

\bibitem{jalu21}
M.~S. {J}akobsen and F.~{L}uef.
\newblock {D}uality of {G}abor frames and {H}eisenberg modules.
\newblock {\em {J}ournal of {N}oncommutative {G}eometry}, 14(4):1445--1500,
  2021.

\bibitem{ka68}
Y.~{K}atznelson.
\newblock {\em {A}n {I}ntroduction to {H}armonic {A}nalysis}.
\newblock {N}ew {Y}ork-{L}ondon-{S}ydney-{T}oronto: {J}ohn {W}iley and {S}ons,
   1968.

\bibitem{ka04-1}
Y.~{K}atznelson.
\newblock {\em {A}n {I}ntroduction to {H}armonic {A}nalysis. 3rd {C}orr. ed.}
\newblock {C}ambridge {U}niv. {P}ress, 2004.

\bibitem{loni23}
E.~{L}orist and Z.~{N}ieraeth.
\newblock {B}anach function spaces done right.
\newblock {\em {I}ndagationes {M}athematicae} 35(2): 247-258, 2024.

\bibitem{lu55}
W.~A.~J. {L}uxemburg.
\newblock {\em {B}anach {F}unction {S}paces}.
\newblock PhD thesis, 1955.

\bibitem{luza71-1}
W.~A.~J. {L}uxemburg and A.~C. {Z}aanen.
\newblock {\em {R}iesz {S}paces {V}ol {I}}.
\newblock {N}orth-{H}olland {P}ublishing {C}ompany. {A}msterdam-{L}ondon, 1971.

\bibitem{pe76}
J.~{P}eetre.
\newblock {\em {N}ew {T}houghts on {B}esov {S}paces}.
\newblock {D}uke {U}niversity {M}athematics {S}eries, {N}o. 1. {M}athematics
  {D}epartment, {D}uke {U}niversity, 1976.

\bibitem{re68}
H.~{R}eiter.
\newblock {\em {C}lassical {H}armonic {A}nalysis and {L}ocally {C}ompact
  {G}roups}.
\newblock {C}larendon {P}ress, {O}xford, 1968.

\bibitem{rest00}
H.~{R}eiter and J.~D. {S}tegeman.
\newblock {\em {C}lassical {H}armonic {A}nalysis and {L}ocally {C}ompact
  {G}roups. 2nd ed.}
\newblock {C}larendon {P}ress, {O}xford, 2000.

\bibitem{re89}
H.~J. {R}eiter.
\newblock {\em {M}etaplectic {G}roups and {S}egal {A}lgebras}.
\newblock {LN} {M}ath.  {S}pringer, {B}erlin, 1989.

\bibitem{sa64-1}
W.~{S}argent.
\newblock {O}n sectionally bounded {B}{K}-spaces.
\newblock {\em Math. Z.}, 83:57--66, 1964.

\bibitem{sa66}
W.~{S}argent.
\newblock {O}n compact matrix transformations between sectionally bounded
  {B}{K}-spaces.
\newblock {\em J. London Math. Soc.}, 41:79--87, 1966.

\bibitem{hosayaya17}
Y.~{S}awano, K.-P. {H}o, D.~{Y}ang, and S.~{Y}ang.
\newblock {H}ardy spaces for ball quasi-{B}anach function spaces.
\newblock {\em {D}issertationes {M}athematicae}, 525:1--102, 2017.

\bibitem{st70}
E.~M. {S}tein.
\newblock {\em {S}ingular {I}ntegrals and {D}ifferentiability {P}roperties of
  {F}unctions}.
\newblock {P}rinceton {U}niversity {P}ress, {P}rinceton, {N}.{J}., 1970.

\bibitem{tayayuzh21}
J.~{T}ao, D.~{Y}ang, W.~{Y}uan, and Y.~{Z}hang.
\newblock {C}ompactness characterizations of commutators on ball {B}anach
  function spaces.
\newblock {\em {P}otential {A}nalysis}, pages 1--35, 2021.

\bibitem{te11}
V.~N. {T}emlyakov.
\newblock {\em {G}reedy {A}pproximation}.
\newblock {C}ambridge {M}onographs on {A}pplied and {C}omputational
  {M}athematics ({N}o. 20). {C}ambridge {U}niversity {P}ress, 2011.

\bibitem{tr83}
H.~{T}riebel.
\newblock {\em {T}heory of {F}unction {S}paces}, vol.~78 {\em {M}onographs
  in {M}athematics}.
\newblock {B}irkh{\"a}user, {B}asel, 1983.

\bibitem{tr92}
H.~{T}riebel.
\newblock {\em {T}heory of {F}unction {S}paces {I}{I}}.
\newblock {M}onographs in {M}athematics 84. {B}irkh{\"a}user, {B}asel, 1992.

\bibitem{tr20}
H.~{T}riebel.
\newblock {\em {T}heory of {F}unction {S}paces {I}{V}}, volume 107 of {\em
  {M}onographs in {M}athematics}.
\newblock {B}irkhauser/{S}pringer, {C}ham, 20.

\bibitem{tr01}
H.~{T}riebel.
\newblock {\em {T}he {S}tructure of {F}unctions.}
\newblock {B}irkh{\"a}user, {B}asel, 2001.

\bibitem{tr16}
H.~{T}riebel.
\newblock {\em {T}empered {H}omogeneous {F}unction {S}paces}.
\newblock {Z}{\"u}rich: {E}uropean {M}athematical {S}ociety ({E}{M}{S}), 2016.

\bibitem{we64}
A.~{W}eil.
\newblock {S}ur certains groupes d\'{ }op{\'e}rateurs unitaires.
\newblock {\em Acta Math.}, 111:143--211, 1964.

\bibitem{za67}
A.~C. {Z}aanen.
\newblock {\em {I}ntegration}.
\newblock {N}orth--{H}olland, {A}msterdam, {C}ompletely revised edition of:
  {A}n {I}ntroduction to the {T}heory of {I}ntegration edition, 1967.

\bibitem{za83-1}
A.~C. {Z}aanen.
\newblock {\em {R}iesz {S}paces {I}{I}}.
\newblock {N}orth-{H}olland {M}athematical {L}ibrary, {V}ol. 30. {A}msterdam -
  {N}ew {Y}ork - {O}xford: {N}orth-{H}olland {P}ubl. {C}ompany. {X}{I},
  1983.

\bibitem{yayuzh23}
C.~{Z}hu, D.~{Y}ang, W.~{Y}uan.
\newblock {G}eneralized {B}rezis--{S}eeger--{V}an {S}chaftingen--{Y}ung
  formulae and their applications in ball {B}anach {S}obolev spaces.
\newblock {\em {C}alc. {V}ar. {P}art.{D}iff. {E}qu.},
  62(8):{P}aper 234, 76, 2023.

\end{thebibliography}

\end{document}